\documentclass[12pt]{article}
\usepackage[dvips]{graphicx}
\usepackage{amssymb,color}
\usepackage[centertags]{amsmath}
\usepackage{amsfonts}

\def\BBox{\kern  -0.2cm\hbox{\vrule width 0.2cm height 0.2cm}}

\title{Combinatorial Space Tiling}
\author{
Egon Schulte
\thanks{Supported by NSF-grant 
DMS--0856675, schulte@neu.edu}\\
{\small Department of Mathematics}\\
{\small Northeastern University, Boston, USA}\\ }

\begin{document}
\maketitle

\begin{abstract}
\noindent
The present article studies combinatorial tilings of Euclidean or spherical spaces by polytopes, serving two main purposes:\  first, to survey some of the main developments in combinatorial space tiling; and second, to highlight some new and some old open problems in this area. 
\end{abstract}


\indent\indent\indent\indent\indent\;\; monotypic tilings



\section{Introduction}
Tilings, or tessellations, have been of interest to artists, scientists and mathematicians throughout history.  In  mathematical investigations, almost all variants of the fundamental question ``{\em How can a given space be tiled by copies of one or more given shapes?\/}" have been studied in some form or another. The underlying spaces have varied from Euclidean, hyperbolic, or spherical spaces, to surfaces or more general manifolds; and the admissible shapes have ranged from well-behaved polyhedral shapes such as convex polytopes, to rather complicated sets with strange properties. The standard assumption in the tiling literature is that a tiling of an ambient space is assembled from congruent copies of shapes from a preassigned set of prototiles. Thus the term {\em copies\/} is usually taken to mean {\em isometric copies\/}.

In combinatorial space tiling, congruence of the tiles is replaced by combinatorial equivalence of the tiles. The study of such tilings was initiated by Danzer in the mid 1970's with the following key problem: \\[.1in]
\indent \vbox{\em Given a convex $d$-polytope in Euclidean $d$-space $\mathbb{E}^d$ does there exist \\
\indent a locally finite face-to-face tiling of $\mathbb{E}^d$ by convex $d$-polytopes each \\
\indent combinatorially equivalent to the given polytope?}\\[.1in]  
(see Danzer, Gr\"unbaum \& Shephard~\cite{dgs}, and \cite{es2})). Here the relaxation of the congruence requirement for the tiles creates many new possibilities for their metrical shapes. In a {\em monotypic\/} tiling, the tiles are convex polytopes that are all of the same combinatorial type, but there may well be infinitely many congruence (or similarity, or affine) classes of tiles. Unless otherwise noted, we still insist on the convexity of the tiles in order to eliminate degenerate possibilities and somewhat limit the degree of freedom for the choice of shapes. In this article we survey some of the main developments in combinatorial space tiling and discuss some new and some old open problems in this area.

The article is organized as follows. In Section~\ref{basicnot} we introduce some basic terminology. Then in Section~\ref{spacfill} we review important results about isometric space-fillers that place the topological results of the subsequent sections into a broader context. Finally, Section~\ref{spacom} is about combinatorial space-fillers and monotypic tilings, discussing known results and posing several challenging open problems.

\section{Basic notions}
\label{basicnot}

For basic terminology on tilings we refer to G\"unbaum \& Shephard~\cite{gs}, Schattschneider \& Senechal~\cite{ss}, and \cite{es1,schuenc}.

A {\em tiling\/}, or {\em tessellation\/}, $\mathcal{T}$ of Euclidean $d$-space $\mathbb{E}^d$ is a countable family of  closed subsets of $\mathbb{E}^d$, the {\em tiles} of $\mathcal{T}$, which cover $\mathbb{E}^d$ without gaps and overlaps; this is to say that the union of all tiles of $\mathcal{T}$ is the entire space, and that any two distinct tiles do not have interior points in common. We generally assume that $\mathcal{T}$ is {\em locally finite\/}, meaning that every compact subset of $\mathbb{E}^d$ meets only finitely many tiles. It is worth pointing out that for many considerations only the real, not the euclidean, structure of the ambient space is required so that we could just as well have used real $d$-space rather than Euclidean $d$-space; in other words, the metric structure on $\mathbb{E}^d$ (distance, angles, etc.) induced by the scalar product will often play a secondary role.

The best behaved tilings are the face-to-face tilings by convex polytopes. A tiling $\mathcal{T}$ of $\mathbb{E}^d$ by convex $d$-polytopes is said to be {\em face-to-face\/} if the intersection of any two tiles is a face of each tile, possibly the empty face. For any such tiling $\mathcal{T}$ of $\mathbb{E}^d$, the $i$-faces of the tiles are also called the {\em $i$-faces\/} of $\mathcal{T}$ ($i = 0,\ldots,d$). The $d$-faces then are the tiles of $\mathcal{T}$, and the $0$- and $1$-faces are the {\em vertices\/} and {\em edges\/} of $\mathcal{T}$, respectively. In particular, the set of all faces of $\mathcal{T}$, ordered by set-theoretic inclusion and suitably extended by a least face and largest face (the empty set and the entire space), is a lattice called the {\em face lattice\/} of $\mathcal{T}$.
Similar terminology applies to more general tilings in which the tiles are {\em topological\/} $d$-polytopes
(homeomorphic images of convex $d$-polytopes) tiling space in a face-to-face manner.

A tiling $\mathcal{T}$ of $\mathbb{E}^d$ by topological $d$-balls is said to be {\em normal\/} if its
tiles are uniformly bounded (meaning that there exist positive parameters $r$ and $R$ such that each tile contains a Euclidean ball of radius $r$, and is contained in a Euclidean ball of radius $R$) and the intersection of every pair of tiles is a connected set (possibly the empty set). The latter condition on the connectedness is trivially satisfied if the tiles are convex $d$-polytopes.

A tiling $\mathcal{T}$ of $\mathbb{E}^d$ is {\em monohedral\/} if all its tiles are congruent to a single set $T$, the ({\em isometric\/}) {\em prototile\/} of $\mathcal{T}$. A convex $d$-polytopes which is the (isometric) prototile of a monohedral tiling of $\mathbb{E}^d$ is often called an ({\em isometric\/}) {\em space-filler\/} of $\mathbb{E}^d$. More generally, a tiling $\mathcal{T}$ of $\mathbb{E}^d$ by convex $d$-polytopes is said to be {\em monotypic\/} if each tile of $\mathcal{T}$ is combinatorially equivalent to a convex $d$-polytope $T$, the {\em combinatorial prototile\/} of $\mathcal{T}$. In a monotypic tiling, there generally are infinitely many different metrical shapes of tiles, but the tiles are all convex and are combinatorially equivalent to a single combinatorial prototile.

\section{Space-fillers}
\label{spacfill}

The classification of (isometric) space-fillers is one of the main open problems in tiling theory (see \cite{gsbull,gs,es1}). The challenge is already evident in the planar case, which is still unsettled. 

It is well-known that a plane-filler must necessarily be a triangle, quadrangle, pentagon or hexagon. The list of convex plane-fillers comprises all triangles and all quadrangles, three kinds of hexagons, and several kinds of pentagons; however, the completeness of the list of pentagonal prototiles has not yet been established. By contrast, combinatorial tiling of the Euclidean plane is easy:\ the plane admits a monotypic face-to-face tiling by convex $p$-gons for each $p \geq 3$ (however, normality can only be achieved for $p\leq 6$). In fact, it is not hard to see that every regular tessellation $\{p,q\}$ of the hyperbolic plane can be realized by a (non-normal) tiling of the Euclidean plane by convex $p$-gons, $q$ meeting at a vertex (see \cite{es3}).

The space-filler problem is most appealing in $3$ dimensions, and here it is widely open. In fact, it is not even known if the number of facets of a $3$-dimensional space-filler is bounded.  Many interesting space-fillers were discovered by crystallographers as Voronoi regions for lattices or other discrete point sets in $\mathbb{E}^3$, including some spectacular examples with as many as 38 facets (see Engel~\cite{eng}). For every dimension $d$, finite upper bounds on the number of facets do exist for space-fillers that admit an isohedral face-to-face tiling of $\mathbb{E}^d$ (see Delone~\cite{delo} and Tarasov~\cite{ta}); in dimension $3$, the best bound known is currently $378$. Recall here that {\em isohedral tilings\/} are monohedral tilings with a tile-transitive symmetry group. The symmetry group of an isohedral tiling must necessarily be among the finitely many crystallographic groups in $\mathbb{E}^d$. The upper bound of \cite{delo,ta} on the number of facets is based on the fact that the index of the translation subgroup of a crystallographic group is uniformly bounded. The existence of an upper bound implies in particular that, for every dimension $d$, there are only finitely many combinatorial types of space-fillers that admit an isohedral face-to-face tiling of $\mathbb{E}^d$.

However, the upper bound results do not apply to general space-fillers. In fact, Hilbert's famous $18^{th}$ Problem, already posed in 1900, asked in part (for $d=3$) whether there exists a $d$-dimensional polyhedral shape that does admit a monohedral, but no isohedral, tiling of $\mathbb{E}^d$.  Such {\em anisohedral\/} space-fillers do in fact exist in every dimension $d\geq 2$ (see \cite{gsbull,gs,es1,ss}). 

It is worth noting that the global property of isohedrality of a tiling can already be detected locally; in fact, the ``Local Theorem for Tilings" says that a face-to-face tiling of $\mathbb{E}^d$ is isohedral if and only if the large enough neighborhoods of tiles satisfy certain conditions (see Dolbilin \& Schattschneider~\cite{dol}). Similar characterizations  also holds for combinatorial tile-transitivity of monotypic tilings, as well as for combinatorial multihedrality of tilings
(see Dolbilin \& Schulte~\cite{dolsch,dolsch1}). Recall here that a tiling is said to be {\em combinatorially multihedral\/} if its combinatorial automorphism group has only finitely many orbits on the tiles.
 
Additional information is available for more narrowly defined classes of tilings or tiles, for example, lattice tilings. In a {\em lattice tiling\/} the translation vectors associated with the tiles form a lattice. If a convex $d$-polytope tiles $\mathbb{E}^d$ by translation, then it also admits (uniquely) a face-to-face lattice tiling of $\mathbb{E}^d$ (see McMullen~\cite{mcm} and Venkov~\cite{ven}); such space-fillers are called {\em parallelohedra\/}. Parallelohedra are the most basic space-fillers. There are only finitely many distinct combinatorial types of parallelohedra in each dimension; however, the exact numbers are known only for $d = 2$, $3$ or $4$, where the numbers are $2$, $5$ or $52$, respectively.

Among the five Platonic solids, only the cube admits a monohedral tiling of ordinary space $\mathbb{E}^3$; there is just one face-to-face tiling by cubes, the {\em regular cubical tessellation\/}. However, there are many tilings of $\mathbb{E}^3$ by cubes that are not face-to-face and in which stacks of cubes are shifted relative to each other. 

In any dimension $d$, there exists a vast variety of tilings by $d$-cubes, some with rather strange properties. It is remarkable that such tilings need not have a pair of cubical tiles with a common $(d-1)$-face; this discovery in Lagarias \& Shor~\cite{lash} disproved a long-standing conjecture, known as Keller's Conjecture, about the existence of such fully adjacent pairs (see also \cite{gsbull}). On the other hand, lattice tilings by $d$-cubes always have pairs of fully adjacent cubes. The class of lattice cube tilings was completely described by Hajos~\cite{haj}, settling a conjecture by Minkowski.

\section{Combinatorial tiling}
\label{spacom}

In combinatorial space tiling, the congruence of the tiles is replaced by combinatorial equivalence while the convexity of the tiles is still maintained. 
\medskip

\noindent
{\bf Combinatorial space-fillers}
\smallskip

It is rather surprising that in ordinary space $\mathbb{E}^3$ {\it every\/} convex polyhedron is a {\em combinatorial space-filler\/}, that is, a combinatorial prototile of a monotypic tiling by convex polytopes (see \cite{es3}). For a generic convex polyhedron in $\mathbb{E}^3$, the corresponding monotypic tiling requires infinitely many congruence classes of tiles and is generated by a rather complicated inductive process in which large star-shaped patches of tiles are extended to even larger star-shaped patches of tiles. 

In general, however, these tilings will not be face-to-face. In fact, there are many convex polyhedra in $\mathbb{E}^3$ which are {\em nontiles\/}, meaning that they do not occur as combinatorial prototiles of monotypic face-to-face tilings of $\mathbb{E}^3$; this even remains true if the tiles in a tiling are allowed to be topological polytopes (see \cite{es2}). The cuboctahedron and the icosidodecahedron are simple examples of nontiles; in spite of their high degree of symmetry, these are not combinatorial prototiles of monotypic {\em face-to-face\/} tilings of $\mathbb{E}^3$. On the other hand, by a result of Gr\"unbaum, Mani \& Shephard~\cite{gms}, every {\em simplicial\/}  (triangle-faced) convex polyhedron in $\mathbb{E}^3$ does occur as the combinatorial prototile of a monotypic face-to-face tiling of $\mathbb{E}^3$; in general, however, infinitely many congruence classes of tiles are required in their construction. 

It is quite remarkable that the icosahedron, as a simplicial convex polyhedron, does admit a monotypic face-to-face tiling of $\mathbb{E}^3$.  It is not known whether or not there exists any monotypic tiling by icosahedra that has only finitely many congruence classes of tiles. For the two other Platonic solids with triangular faces, a single congruence class of tiles suffices. 

There are several kinds of tilings of $\mathbb{E}^3$ by congruent tetrahedra. For example, the chamber complex of the regular cubical tessellation of $\mathbb{E}^3$ (consisting of the barycentric subdivisions of the cubical tiles) is a face-to-face tiling of $\mathbb{E}^3$ by congruent tetrahedra in which any tetrahedral tile is a fundamental region (fundamental tetrahedron) for the action of the symmetry group of the cubical tessellation on $\mathbb{E}^3$. The question which tetrahedra tile space $\mathbb{E}^3$ has a long history, going back (at least) to Aristotle, who incorrectly claimed that the regular tetrahedron could tile $\mathbb{E}^3$ (see Senechal~\cite{sentet} and Debrunner~\cite{debsim}). The full classification of the tetrahedra that admit a face-to-face tiling by properly congruent tiles was recently obtained in Edmonds~\cite{edm}, thereby establishing that the list of Sommerville~\cite{som} is complete.

From the regular cubical tessellation we can also construct a nice face-to-face tiling of $\mathbb{E}^3$ by congruent (non-regular) octahedra; each pair of adjacent cubes determines an octahedral tile obtained by joining the center of the cubes to their common square face. For recent work on the classification of (topologically tile transitive) tilings of $\mathbb{E}^3$ by topological tetrahedra, octahedra, or cubes, the reader is also referred to Delgado-Friedrichs \& Huson~\cite{dfh}. For applications of tilings by topological Platonic solids or other highly symmetric tiles to the classification of three-periodic nets in crystal chemistry see also Delgado-Friedrichs \& al.~\cite{dal}.

In the next subsection we show that the dodecahedron, although non-simplicial, also admits monotypic face-to-face tilings; in fact, there are such tilings with relatively few congruence classes of tiles.

We mention in passing that there is also an interesting concept of combinatorial aperiodicity that generalizes the well-known concept of geometric aperiodicity of a prototile set (see \cite{esap}).  A finite prototile set of convex polytopes is called {\em combinatorially aperiodic\/} if it admits a locally finite face-to-face tiling by combinatorially equivalent convex copies of the prototiles, but if no such tiling has a combinatorial automorphism of infinite order. 
There are no combinatorially aperiodic prototile set in the plane. It is an open problem whether there are any combinatorially aperiodic prototile sets in dimensions $d\geq 3$ (see \cite{esap}). Even the case of a single prototile has not been settled.
\medskip

\noindent
{\bf Projection techniques}
\smallskip

Before proceeding, it is important to point out a close connection between combinatorial tiling properties of convex $d$-polytopes in $\mathbb{E}^d$ and their appearance as facet types of equifacetted convex $(d+1)$-polytopes in $\mathbb{E}^{d+1}$. Recall here that a convex polytope is {\em equifacetted\/} if all its facets are combinatorially equivalent. Clearly, convex $d$-polytopes that admit face-to-face tilings (with finitely many congruence classes of tiles, or not) do not in general occur as facet types of equifacetted convex $(d+1)$-polytopes.  For example, the regular hexagon tiles the plane, but there are no convex polyhedra in $\mathbb{E}^3$ with only hexagonal faces. However, an equifacetted polytope always produces a tiling, as can be seen by applying projection techniques that generate an infinite sequence of nested Schlegel diagrams converging in the limit to a tiling of the entire space. In fact, every convex $d$-polytope occurring as the facet type of an equifacetted convex $(d+1)$-polytope is also the combinatorial prototile of a monotypic face-to-face tiling of $\mathbb{E}^{d+1}$ (see \cite{es2}). However, it is not known yet if in general there is also a normal such tiling. On the other hand, it may be conjectured that only finitely many congruence classes of tiles are needed, with their number depending on the given prototile.

Again the case of the icosahedron is particularly interesting. It is a long-standing open problem to decide whether or not the icosahedron occurs as the facet type of an equifacetted convex $4$-polytope (see Perles \& Shephard \cite{persh}, Kalai~\cite{kalai1, kalai2} and Barnette~\cite{bar}). On the other hand, as we saw earlier, the icosahedron does admit a monotypic face-to-face tiling of $\mathbb{E}^3$.  All other Platonic solids occur as facet types of equifacetted convex $4$-polytopes. In fact, the tetrahedron, cube, octahedron, and dodecahedron even occur as facets of a convex regular $4$-polytope, namely the $4$-simplex (or $4$-crosspolytope or $600$-cell), the $4$-cube, the $24$-cell, or the $120$-cell, respectively. 

It is known that the truncated icosahedron, also known as the {\em football\/}, is not the facet type of an equifacetted topological $4$-polytope (see Pasini~\cite{pas}). However, there does exist a monotypic face-to-face tiling of $\mathbb{E}^3$ by topological copies of truncated icosahedra; in fact, the truncated icosahedron does even admit a tiling by convex polytopes in hyperbolic $3$-space $\mathbb{H}^3$. The football is an example of a {\em fullerene\/}, that is, a polyhedron with only $3$-valent vertices and only pentagonal or hexagonal faces. See also Dutour Sikiri\'c, Delgado-Friedrichs \& Deza~\cite{ddd} for tilings of $\mathbb{E}^3$ 
by several standard fullerenes.

There is one particular case, which we briefly review here, when certain projection methods combined with reflection group techniques allow some very elegant constructions of monotypic tilings with finitely many congruence classes of tiles. Suppose that a convex $d$-polytope $P$ is realized as the facet type of an equifacetted convex $(d+1)$-polytope $Q$ which has least one $(d+1)$-valent vertex. Then $P$ is the combinatorial prototile of a periodic monotypic face-to-face tiling of $\mathbb{E}^d$ with at most $m-d-1$ congruence classes of tiles, where $m$ is the number of facets of $Q$ (see \cite{es2}). Recall here that a tiling is {\em periodic\/} if its symmetry group has a translation subgroup of rank $3$ generated by three translations in independent directions. This tiling of $\mathbb{E}^d$ can be constructed as follows.

Suppose that $x$ is a $(d+1)$-valent vertex of $Q$, and that $x_{1},\ldots,x_{d+1}$ are the vertices of $Q$ adjacent to $x$. Then $x_{1},\ldots,x_{d+1}$ span a $d$-simplex $T$ whose affine hull strictly separates $x$ from all the other vertices of $Q$. Now project the boundary complex of $Q$, radially from $x$, onto the affine hull of $T$. This yields a face-to-face ``tiling" $\mathcal{C}$ of $T$ by convex $d$-polytopes each isomorphic to $P$. In fact, the $m-d-1$ facets of $Q$ that do not contain $x$ project onto the tiles of $\mathcal{C}$, and those that do contain $x$ become the $(d-1)$-faces of $T$. Thus $\mathcal{C}$ is the image of the anti-star of $x$ in the boundary complex of $Q$. 

Next we exploit the fact, pointed out above for $d=3$, that each $d$-simplex in the chamber complex of the regular cubical tessellation of $\mathbb{E}^d$ is a fundamental region for the action of the symmetry group of the tessellation on $\mathbb{E}^d$. Suppose $T'$ is any such $d$-simplex. Now map $T$ affinely onto $T'$ to generate a ``tiling" $\mathcal{C}'$ of $T'$ from the tiling $\mathcal{C}$ of $T$. Then the tiles in $\mathcal{C}'$ are affine images of those in $\mathcal{C}$ and hence are convex polytopes isomorphic to $P$. Finally, then, a face-to-face tiling of the entire space $\mathbb{E}^d$ is obtained from $\mathcal{C}'$ by applying all symmetries of the cubical tessellation. The various copies of $\mathcal{C}'$ under these symmetries fit together in a face-to-face manner; this follows from basic properties of the action of the symmetry group, viewed as a reflection (Coxeter) group generated by the reflections in the walls of $T'$, on the chamber complex. The number of congruence classes of tiles in the resulting space tiling is just that of $\mathcal{C}'$, and hence is at most $m-d-1$.

The most prominent example to which the previous construction applies is the regular $120$-cell $Q=\{5,3,3\}$ in $\mathbb{E}^4$, which has $120\,(=m)$ dodecahedral facets and $600$ vertices, each with tetrahedral vertex-figure (see Coxeter~\cite{crp}). Due to the high degree of symmetry of $Q$, the corresponding tetrahedron $T$ is regular and $\mathcal{C}$ actually has considerably fewer congruence classes of tiles than $116 \,(=m-d-1)$. However, in applying the affine transformation of $T$ onto $T'$ much of the gain in the number of congruence classes is lost again; this is due to the fact that, while $T$ has maximum possible symmetry, $T'$ has only one non-trivial symmetry (a half-turn through a suitable pair of opposite edges). The author has not established the precise number of congruence classes of tiles in the resulting space  tiling, but it is clear that the exact number is considerably smaller than $116$. 
\medskip

\noindent
{\bf Exploiting Coxeter group methods}
\smallskip

In trying to significantly decrease the number of congruence classes of dodecahedral tiles, a good strategy consists of replacing the chamber complex of the cubical tessellation of $\mathbb{E}^3$ by that of another affine Coxeter group such as the Coxeter group $\widetilde{A}_3$ (with a circular diagram on four nodes and with unmarked branches); the presence of additional diagram symmetries translates into more Euclidean symmetries of the corresponding fundamental chamber and hence in less loss in the transition from $\mathcal{C}$ to $\mathcal{C}'$. Employing $\widetilde{A}_3$ should further reduce the number of congruence classes of dodecahedral tiles. However, it is clear from the onset that there are limits to any such reduction and that no tilings with a very small number of congruence classes of tiles can be found this way. It is an interesting open problem to determine the minimum number of congruence classes of tiles that a monotypic (face-to-face) tiling by dodecahedra can have? In principle, this number of classes could be as low as $2$ or $3$, and could even be $1$?

In this context it may be helpful to remember that the dodecahedron is the tile in two regular tessellations of hyperbolic $3$-space $\mathbb{H}^3$, namely in $\{5,3,4\}$ (with octahedral vertex-figures) and in$\{5,3,5\}$ (with icosahedral vertex-figures); see \cite{coxhyp} for more details. Topologically, these give us combinatorially regular tessellations of Euclidean $3$-space $\mathbb{E}^3$ with topological dodecahedra. It would be interesting to know if the entire space $\mathbb{E}^3$ can be tiled in a locally finite face-to-face manner by convex dodecahedra, and if one can further achieve affine equivalence of the tiles. The following observation hints that it may indeed be possible to achieve convexity of the tiles. 

In fact, the two regular tessellations of $\mathbb{H}^3$ with dodecahedral tiles can be realized by isomorphic ``tilings" of a $3$-dimensional open unit ball by convex dodecahedra. A proof of this fact exploits the canonical representation of the Coxeter (symmetry) groups of $\{5,3,4\}$ and $\{5,3,5\}$ in real $4$-space, in particular the action of these groups as linear groups on the corresponding chamber complex (see McMullen \& Schulte~\cite[Section 3D]{arp}). Here this chamber complex consists of $4$-dimensional simplicial cones tessellating a $4$-dimensional convex cone, known as the Tits cone; this cone can be viewed as hyperbolic $3$-space, represented in the projective model. The chambers (cones) can be grouped together in sets of $120$, each determining a $4$-dimensional convex dodecahedral cone, and these dodecahedral cones fit together to form the tiles in a tiling of the Tits cone that is combinatorially equivalent to $\{5,3,4\}$ or $\{5,3,5\}$. (This is the geometric model of the universal regular $4$-polytope $\{5,3,4\}$ or $\{5,3,5\}$, respectively, described in \cite[Theorem 3D7]{arp}.)  When the Tits cone is cut by a suitable hyperplane in $4$-space, a tiling of an open $3$-ball by convex dodecahedra is obtained. After rescaling, this gives a tiling of the open unit $3$-ball combinatorially equivalent to $\{5,3,4\}$ or $\{5,3,5\}$.

Similar considerations apply to other regular tessellations in hyperbolic $3$-space or $4$-space. For example, for  the regular tessellation $\{3,5,3\}$ of $\mathbb{H}^3$ (with icosahedral tiles and dodecahedral vertex-figures) we arrive at tilings of an open unit $3$-ball by convex icosahedra. Similarly, from the three regular tessellations $\{5,3,3,3\}$, $\{5,3,3,4\}$, and $\{5,3,3,5\}$ of $\mathbb{H}^4$ (with $120$-cells as tiles, and with $4$-simplices, $4$-crosspolytopes, or $600$-cells as vertex-figures, respectively) we obtain tessellations of the open unit $4$-ball by convex $120$-cells.

\noindent
{\bf More on nontiles}
\smallskip

Although the degree of freedom in designing monotypic tilings is much larger than for monohedral tilings, there still are major obstructions arising from the combinatorics of the prototiles or the convexity requirement for the tiles. In particular, for $d\geq 3$, there are many different kinds of convex $d$-polytopes that are nontiles and hence also do not occur as facet types of equifacetted convex $(d+1)$-polytopes (see \cite{es2}). Many examples are nontiles in a strong, topological sense, in that they even do not admit monotypic face-to-face tilings by topological polytopes. Other nontile criteria exploit the convexity of the tiles in a tiling and lead to nontiles in the original sense. 

It is quite spectacular that, in spite of its high degree of symmetry, the $d$-crosspolytope is a nontile when $d\geq 7$. This is in sharp contrast to the fact that its dual, the $d$-cube, is the simplest tile imaginable and gives rise to a regular tessellation in any dimension $d$. In dimension $4$, the crosspolytope even admits a regular tessellation of space, as does the regular $24$-cell.  The combinatorial tiling properties of the $120$-cell and $600$-cell, as well as those of the $6$-crosspolytope, do not seem to be known (see also \cite{es2}); on the other hand, the $5$-crosspolytope does not tile $\mathbb{E}^5$ in a face-to-face manner (see Kalai~\cite{kalai1}). 

On the other hand, many basic figures in $\mathbb{E}^3$ such as $p$-gonal prisms, pyramids, or bipyramids all admit monotypic face-to-face tilings (see \cite{es2}). For example, nice monotypic face-to-face tilings of $\mathbb{E}^3$ with only finitely many congruence classes of $p$-gonal prisms can be derived by the projection techniques applied to the (simple) cartesian product of two regular $p$-gons in complementary planes of $\mathbb{E}^4$.

The tiling properties of simple convex polytopes are not well understood. Recall here that a convex $d$-polytope is called {\em simple\/} if all its vertices are $d$-valent. In particular, no simple nontiles are known in any dimension. However, there are simple convex polytopes in dimensions $3$ and $4$ which are {\em nonfacets\/}, meaning that they are not the facet types of equifacetted convex polytopes (of dimensions $4$ or $5$, respectively); see Perles-Shephard~\cite{persh} and Barnette~\cite{bar}. It may be conjectured that every simple convex $3$-polytope is the prototile of a monotypic face-to-face tiling of $\mathbb{E}^3$ by topological polytopes. However, the author expects the answer to the corresponding question in higher dimensions to be negative.

In dealing with tilings by topological polytopes, the face-to-face condition often eliminates pathological situations. In fact, as the following construction shows, every convex $d$-polytope is the prototile of a highly degenerate monotypic tiling of $\mathbb{E}^d$ by topological polytopes that is not fact-to-face in general. First recall that the boundary complex of any convex $d$-polytope is isomorphic to a refinement of the boundary complex of the $d$-simplex; see Gr\"unbaum~\cite{gru} for details. In dimension $3$ this says that, up to homeomorphism, the boundary complex of the polytope can viewed as decomposing the boundary of a tetrahedron in a suitable way. When this tetrahedron, along with its boundary suitably decomposed by the polytope boundary complex, is taken to be the fundamental tetrahedron in the chamber complex of the regular cubical tessellation in $\mathbb{E}^3$, and all symmetries of the cubical tessellation are applied, we obtain a monotypic tiling of $\mathbb{E}^3$ by topological polytopes in which every tile is isomorphic to the given polytope. However, the face-to-face property fails in a major way, unless the original polytope itself was a tetrahedron. These arguments extend to any dimension.
\medskip

\noindent
{\bf Tilings of $3$-space by handlebodies}
\smallskip

A largely unexplored direction in combinatorial tiling is the study of tilings by handlebodies in Euclidean $3$-space $\mathbb{E}^3$ or the unit $3$-sphere $\mathbb{S}^3$. In our previous discussion the tiles were always assumed to be convex (or topological) $d$-polytopes and hence to be topologically spherical (or rather, homeomorphic to $d$-balls). Interesting new possibilities arise if handlebodies of genus $g$ are permitted as tiles. Examples of such tilings have been described in, for example, Adams~\cite{adams,adams1}, Coxeter \& Shephard~\cite{coxshe}, Debrunner~\cite{debr}, Kuperberg~\cite{kup}, and \cite{schusp}.

From the combinatorial perspective we are primarily interested in tiles whose bounding (orientable) surface of genus $g$ carries a $2$-dimensional boundary complex that is an abstract polyhedron in the sense of McMullen \& Schulte~\cite{arp} (or an orientable map in the sense of Coxeter \& Moser~\cite{cm}). The  handlebodies in a tiling, each equipped with its boundary complex, then are required to fit together in a face-to-face manner to tile $\mathbb{E}^3$, again with local finiteness understood. Each such face-to-face tiling of $\mathbb{E}^3$ or $\mathbb{S}^3$ by handlebodies has itself a natural structure of a complex in which the $3$-faces correspond to the tiles. Frequently this complex is an abstract polytope of rank $4$, and the tiling is a topological model for this polytope (see \cite[Section 6B]{arp}). A more precise formulation of topological modeling of abstract polytopes, as well as a number of interesting examples (including tessellation of more general topological spaces), can be found Brehm, K\"uhnel \& Schulte~\cite{bks}.

For example, the unit $3$-sphere $\mathbb{S}^3$ can be tessellated by $20$ toroidal (picture frame-like) handlebodies, each equipped with a boundary complex consisting of $9$ squares arranged in a $3\times 3$ fashion. This tiling, independently discovered by Gr\"unbaum~\cite{polstra} and Coxeter \& Shephard~\cite{coxshe}, is a $3$-dimensional spherical model for the universal locally toroidal regular $4$-polytope $\{\{4,4\}_{(3,0)},\{4,3\}\}$. This polytope has $20$ toroidal facets $\{4,4\}_{(3,0)}$, $30$ vertices (with cubes $\{4,3\}$ as vertex-figures), and combinatorial automorphism group $S_{6}\times C_2$ (see \cite[Section 10B]{arp}). 

It would be very interesting to know which of the (universal) locally toroidal regular $4$-polytopes admit tilings of $\mathbb{E}^3$ or $\mathbb{S}^3$ as topological models? For example, can every universal regular $4$-polytope $\{\{4,4\}_{(m,0)},\{4,3\}\}$, with $m\geq 4$, be realized by a tiling of $\mathbb{E}^3$ with toroidal handlebodies?  
This particular question concerns polytopes of type $\{4,4,3\}$, but there also similar such questions for the types $\{6,3,3\}$, $\{6,3,4\}$ and $\{6,3,5\}$ (with hexagon-faced toroidal handlebodies). Furthermore, relating compactness of an ambient space to the finiteness of a polytope, it can be asked if every finite universal locally toroidal regular $4$-polytope (with spherical vertex-figures) can be modeled by a tiling on $\mathbb{S}^3$ with toroidal handlebodies? Or, even more generally, which finite locally toroidal abstract regular $4$-polytopes can be modeled by a tiling on $\mathbb{S}^3$? These basic questions can also be asked for handlebodies of higher genus. It seems that no general results in this direction are known.

\nopagebreak

 \end{document}